\documentclass[12pt]{article}
\usepackage{amsfonts,amsmath,amsxtra}
\usepackage{latexsym}
\usepackage{amssymb}

\usepackage{comment}

\def\hybrid{\topmargin 0pt      \oddsidemargin 0pt
        \headheight 0pt \headsep 0pt
        \textwidth 16.5cm
        \textheight 23cm
        \marginparwidth 0.0in
        \parskip 5pt plus 1pt   \jot = 1.5ex}
\catcode`\@=11
\def\marginnote#1{}
\newcount\hour
\newcount\minute
\newtoks\amorpm
\hour=\time\divide\hour by60 \minute=\time{\multiply\hour by60
\global\advance\minute by-\hour}
\edef\standardtime{{\ifnum\hour<12 \global\amorpm={am}%
        \else\global\amorpm={pm}\advance\hour by-12 \fi
        \ifnum\hour=0 \hour=12 \fi
      \number\hour:\ifnum\minute<10 0\fi\number\minute\the\amorpm}}
\edef\militarytime{\number\hour:\ifnum\minute<10 0\fi\number\minute}

\def\draftlabel#1{{\@bsphack\if@filesw {\let\thepage\relax
   \xdef\@gtempa{\write\@auxout{\string
      \newlabel{#1}{{\@currentlabel}{\thepage}}}}}\@gtempa
   \if@nobreak \ifvmode\nobreak\fi\fi\fi\@esphack}
        \gdef\@eqnlabel{#1}}
\def\@eqnlabel{}
\def\@vacuum{}
\def\draftmarginnote#1{\marginpar{\raggedright\scriptsize\tt#1}}

\def\draft{\oddsidemargin -0.1truein
        \def\@oddfoot{\sl preliminary draft \hfil
        \rm\thepage\hfil\sl\today\quad\militarytime}
        \let\@evenfoot\@oddfoot \overfullrule 3pt
        \let\label=\draftlabel
        \let\marginnote=\draftmarginnote
\def\@eqnnum{{\rm (\theequation)}
\rlap{\kern\marginparsep\tt\@eqnlabel}%
\global\let\@eqnlabel\@vacuum}  }
\newfont{\Bbbb}{msbm7 scaled 1\@ptsize00}
\newcommand{\zs}{\raise-1pt\hbox{$\mbox{\Bbbb Z}$}}

scaled 1\@ptsize00
\font\sevenmsa=msam6 %scaled 1\@ptsize00
\newfam\msafam
\textfont\msafam=\sevenmsa
\def\hexnumber@#1{\ifnum#1<10 \number#1\else
\ifnum#1=10 A\else\ifnum#1=11 B\else\ifnum#1=12 C\else \ifnum#1=13
D\else\ifnum#1=14 E\else\ifnum#1=15 F\fi\fi\fi\fi\fi\fi\fi}
\def\msa@{\hexnumber@\msafam}
\def\llcorner{\delimiter"4\msa@78\msa@78 }
\def\lrcorner{\delimiter"5\msa@79\msa@79 }
\mathchardef\blacktriangleright="3\msa@49
\mathchardef\blacktriangleleft="3\msa@4A \font\tenmsb=msbm10 scaled
1\@ptsize00
\newfam\msbfam
\textfont\msbfam=\tenmsb \scriptfont\msbfam=\tenmsb
\newdimen\Squaresize \Squaresize=14pt
\newdimen\Thickness \Thickness=0.5pt

\def\Square#1{\hbox{\vrule width \Thickness
   \vbox to \Squaresize{\hrule height \Thickness\vss
      \hbox to \Squaresize{\hss#1\hss}
   \vss\hrule height\Thickness}
\unskip\vrule width \Thickness} \kern-\Thickness}

\def\Vsquare#1{\vbox{\Square{$#1$}}\kern-\Thickness}

\def\numberbysection{\@addtoreset{equation}{section}
        \def\theequation{\thesection.\arabic{equation}}}
\numberbysection

\renewcommand{\theequation}{\thesection.\arabic{equation}}
\def\titlepage{\@restonecolfalse\if@twocolumn\@restonecoltrue\onecolumn
     \else \newpage \fi \thispagestyle{empty}\c@page\z@
        \def\thefootnote{\fnsymbol{footnote}} }

\def\endtitlepage{\if@restonecol\twocolumn \else  \fi
        \def\thefootnote{\arabic{footnote}}
        \setcounter{footnote}{0}}  %\c@footnote\z@ }
\relax

\hybrid
\makeatletter
\newdimen\normalarrayskip            % skip between lines
\newdimen\minarrayskip               % minimal skip between lines
\normalarrayskip\baselineskip \minarrayskip\jot
\newif\ifold             \oldtrue            \def\new{\oldfalse}
\def\arraymode{\ifold\relax\else\displaystyle\fi}%mode of array enrties
\def\eqnumphantom{\phantom{(\theequation)}} % ight phantom in eqnarray
\def\@arrayskip{\ifold\baselineskip\z@\lineskip\z@
     \else
     \baselineskip\minarrayskip\lineskip1\baselineskip\fi}

\def\@arrayclassz{\ifcase \@lastchclass \@acolampacol \or
\@ampacol \or \or \or \@addamp \or
   \@acolampacol \or \@firstampfalse \@acol \fi
\edef\@preamble{\@preamble
  \ifcase \@chnum
     \hfil$\relax\arraymode\@sharp$\hfil
     \or $\relax\arraymode\@sharp$\hfil
     \or \hfil$\relax\arraymode\@sharp$\fi}}

\def\@array[#1]#2{\setbox\@arstrutbox=\hbox{\vrule
     height\arraystretch \ht\strutbox
     depth\arraystretch \dp\strutbox
width\z@}\@mkpream{#2}\edef\@preamble{\halign \noexpand\@halignto
\bgroup \tabskip\z@ \@arstrut \@preamble \tabskip\z@ \cr}%
\let\@startpbox\@@startpbox \let\@endpbox\@@endpbox
  \if #1t\vtop \else \if#1b\vbox \else \vcenter \fi\fi
  \bgroup \let\par\relax
  \let\@sharp##\let\protect\relax
  \@arrayskip\@preamble}
\def\eqnarray{\stepcounter{equation}%
              \let\@currentlabel=\theequation
              \global\@eqnswtrue
              \global\@eqcnt\z@
              \tabskip\@centering              %formulae  centering
              \let\\=\@eqncr
              $$%
            \halign to \displaywidth  \bgroup
             \eqnumphantom \@eqnsel
      \hskip\@centering                               %right tab%
    $\displaystyle  \tabskip\z@ {##}$%
    &\global\@eqcnt\@ne \hskip 2\arraycolsep
         $ \displaystyle  \arraymode{##}$\hfil
    &\global\@eqcnt\tw@ \hskip 2\arraycolsep
         $\displaystyle\tabskip\z@{##}$\hfil
         \tabskip\@centering
    &{##}\tabskip\z@\cr}
\makeatother
\newtheorem{te}{Theorem}[section]%Usage:\begin{te}Statement\end{te}

\newtheorem{prop}{Proposition}[section]           %  ETC ...

\newtheorem{lem}{Lemma}[section]

\newcommand\bqa{\begin{eqnarray}}
\newcommand\eqa{\end{eqnarray}}
\def\be{\begin{eqnarray}\new\begin{array}{cc}}
\def\ee{\end{array}\end{eqnarray}}

\def\beq{\begin{equation}}
\def\eeq{\end{equation}}
\def\bse{\begin{subequations}}                %%%SUBEQUATIONS
\def\ese{\end{subequations}}
\def\bp{\begin{pmatrix}}
\def\ep{\end{pmatrix}}

\def\proof{\noindent {\it Proof}. }
\def\stack#1#2{\raise0.7pt\hbox{$\mathrel{\mathop{#2}\limits^{#1}}$}}
\def\tr{\triangleright}
\def\tl{\triangleleft}
\def\sem{\mathsurround=0pt \raise1pt
\hbox{$\scriptscriptstyle>\!\!$}\:\!\!\tl}
\def\mes{\mathsurround=0pt \tr\!\:\!\raise0.8pt
\hbox{$\scriptscriptstyle\!\!<$}\,}
\def\]{\mathsurround=0pt ]\raise-2pt\hbox{$_\ast$}}
\def\<{\langle}
\def\>{\rangle}

\def\we{\raise-1pt\hbox{$\,\stackrel{\wedge}{,}\,$}}
\def\tr{{\rm tr}\,}

\newcounter{pac}[section]

\newcounter{pacc}[subsection]

0
\begin{document}

%%%%%%%%%%%% Title page %%%%%%%%%%%%%%%%%%%%%
\setcounter{pac}{0}
\setcounter{footnote}0

\begin{center}

\phantom.
\bigskip%\bigskip\bigskip\bigskip\bigskip\bigskip
%{\hfill{\normalsize ITEP-TH-30/04}\\
%\hfill{\normalsize TCD-MATH-04-15}\\
%\hfill{\normalsize HMI-04-04}\\
%\hfill{\normalsize MPIM-....}\\
%[15mm]\Large\bf

{\Large\bf On explicit realization of algebra of complex divided powers of $U_{q}(\mathfrak{sl}(2))$}

\vspace{1cm}

\bigskip\bigskip

{\large  Pavel Sultanich \footnote {E-mail:  sultanichp@gmail.com}},\\
\bigskip
{\it Moscow Center for Continuous Mathematical Education, 119002, Bolshoy Vlasyevsky Pereulok 11, Moscow, Russia
}\\
\bigskip

\end{center}

%\date{}

%\maketitle

%\renewcommand{\abstractname}{}

\begin{abstract}
%\centering{\noindent {\bf }}
\noindent

In this note we prove that the explicit realization of arbitrary complex powers of generators of quantum group $U_{q}(\mathfrak{sl}(2))$ satisfies all the commutation relations of the algebra of complex powers, including the generalized Kac's identity which was announced in our previous paper. It turns out that the latter identity in this realization is equivalent to $6-9$ integral identity on quantum dilogarithm.

\end{abstract}
\vspace{5 mm}

\section{Introduction}

Let $\mathfrak{g}$ be a Lie algebra, let $q = e^{\pi\imath b^{2}}$, $b^{2}\in (0;1)$, and $b^{2}$ is irrational. To any such Lie algebra one can associate a Hopf algebra $U_{q}(\mathfrak{g})$ called quantum group \cite{ChPr}.
The notion of the modular double of quantum group was introduced by Faddeev \cite{F2} in the case of $U_{q}(\mathfrak{sl}(2))$ as a tensor product of two quantum groups $U_{q}(\mathfrak{sl}(2))$ and its modular dual $U_{\tilde{q}}(\mathfrak{sl}(2))$, $\tilde{q} = e^{\pi\imath b^{-2}}$. The modular double appears in many areas of mathematical physics such as Liouville theory \cite{PT}, \cite{FKV}, relativistic Toda model \cite{KLSTS} and others. Further progress in the study of modular double has been made in papers \cite{BT}, \cite{PT1}, where the special class of representations of $U_{q}(\mathfrak{sl}(2))$ has been considered. These representations exhibit a duality under the exchange $b \leftrightarrow b^{-1}$ and they are simultaneously representations of the modular dual group $U_{\tilde{q}}(\mathfrak{sl}(2))$. Furthermore, the generators of the dual group are related to the generators of the original via the so-called transcendental relations. Specifically, let $K$, $E$, $F$ be the usual generators of $U_{q}(\mathfrak{sl}(2))$ and  $\tilde{K}$, $\tilde{E}$, $\tilde{F}$ be the generators of the dual group $U_{\tilde{q}}(\mathfrak{sl}(2))$. Define the rescaled versions of some of these generators
\begin{equation}\label{rescaled E}
    \mathcal{E} = -\imath (q-q^{-1})E,
\end{equation}
\begin{equation}\label{rescaled F}
    \mathcal{F} = -\imath (q-q^{-1})F,
\end{equation}
and the dual ones
\begin{equation}\label{rescaled E}
    \tilde{\mathcal{E}} = -\imath (\tilde{q}-\tilde{q}^{-1})\tilde{E},
\end{equation}
\begin{equation}\label{rescaled F}
    \tilde{\mathcal{F}} = -\imath (\tilde{q}-\tilde{q}^{-1})\tilde{F}.
\end{equation}
Then these rescaled generators are realized by positive self-adjoint operators and satisfy the transcendental relations \cite{BT}
\begin{equation}
    \tilde{K} = K^{b^{-2}},
\end{equation}
\begin{equation}
    \tilde{\mathcal{E}} = \mathcal{E}^{b^{-2}},
\end{equation}
\begin{equation}
    \tilde{\mathcal{F}} = \mathcal{F}^{b^{-2}}.
\end{equation}

After identification of the generators of the dual quantum group with certain powers of generators of the original, one is naturally led to consider arbitrary complex powers of generators which will form a basis of the bigger Hopf algebra.

In \cite{Su} we presented  complete set of defining relations\footnote{Part of them was obtained previously in e.g. \cite{Ip}.} of the Hopf algebra of arbitrary complex powers in the case of simply-laced Lie algebras.

In this note we study explicit realization of the Hopf algebra of arbitrary complex powers of generators of quantum group $U_{q}(\mathfrak{sl}(2,\mathbb{R}))$ via positive self-adjoint operators. We check that all the defining relations, including the generalized Kac's identity hold in this realization.

The plan of the paper is as follows. In Section 2 we recall the definition of quantum group $U_{q}(\mathfrak{sl}(2))$, special function $G_{b}(x)$, called quantum dilogarithm and its properties which plays an important role in the theory of positive principal series representations. In Section 3 the Hopf algebra of arbitrary complex analogs of divided powers of generators has been constructed, which is the generalization of the algebra of divided powers $\frac{X^{n}}{[n]_{q}!}$ studied by Lusztig (see e.g. \cite{Lu 1}).  We study explicit realization of the Hopf algebra generated by divided complex powers of generators of $U_{q}(\mathfrak{sl}(2))$. In the Theorem 3.1 we summarize all defining relations in the Hopf algebra of complex powers of $U_{q}(\mathfrak{sl}(2))$.

{\bf Acknowledgements:} The research was supported by    RSF  (project 16-11-10075).

\newpage

\section{Preliminaries}

Let us start by recalling the definition of a quantum group $U_{q}(\mathfrak{sl}(2))$ \cite{ChPr}.
$U_{q}(\mathfrak{sl}(2))$ $(q = e^{\pi\imath b^{2}}$, $b^{2}\in \mathbb{R}\setminus \mathbb{Q})$ is a Hopf algebra with generators $E$, $F$, $K = q^{H}$ and relations
\begin{equation}
    KK^{-1} = K^{-1}K,
\end{equation}
\begin{equation}
    KE = q^{2}EK,
\end{equation}
\begin{equation}
    KF = q^{-2}FK,
\end{equation}
\begin{equation}
    EF - FE = \frac{K - K^{-1}}{q-q^{-1}}.
\end{equation}
Coproduct is given by
\begin{equation}
\Delta E = E\otimes 1 + K^{-1}\otimes E,
\end{equation}
\begin{equation}
\Delta F = 1\otimes F + F\otimes K,
\end{equation}
\begin{equation}
\Delta K = K\otimes K.
\end{equation}

Non-compact quantum dilogarithm $G_{b}(z)$ is a special function introduced in \cite{F1} (see also \cite{F0}, \cite{FKV}, \cite{V}, \cite{Ka}, \cite{KLSTS}, \cite{BT}). It is defined as follows
\begin{equation}
\log G_{b}(z) = \log\bar{\zeta}_{b} - \int\limits_{\mathbb{R}+\imath 0} \frac{dt}{t}\frac{e^{zt}}{(1-e^{bt})(1-e^{b^{-1}t})},
\end{equation}
where $Q = b+b^{-1}$ and $\zeta_{b} = e^{\frac{\pi\imath}{4} + \frac{\pi\imath(b^{2}+b^{-2})}{12}}$. Note, that $G_{b}(z)$ is closely related to the double sine function $S_{2}(z|\omega_{1},\omega_{2})$, see eq.(A.22) in \cite{KLSTS}.

Below we outline some properties of $G_{b}(z)$, for details see appendix.\\*
1. The function $G_{b}(z)$ has simple poles and zeros at the points
\begin{equation}
    z = -n_{1}b -n_{2}b^{-1},
\end{equation}
\begin{equation}
    z = Q +n_{1}b + n_{2}b^{-1},
\end{equation}
respectively, where $n_{1}$,$n_{2}$ are nonnegative integer numbers.\\*
2. $G_{b}(z)$ has the following asymptotic behavior:
\begin{equation}
 G_{b}(z) \sim
 \begin{cases} \bar{\zeta}_{b}, Im z \rightarrow +\infty ,\\ \zeta_{b} e^{\pi\imath z(z-Q)}, Im z \rightarrow -\infty . \end{cases}
\end{equation}
3. Functional equation:
\begin{equation}
G_{b}(z +b^{\pm 1}) = (1-e^{2\pi\imath b^{\pm 1}z})G_{b}(z).
\end{equation}
4. Reflection formula:
\begin{equation}
G_{b}(z)G_{b}(Q-z) = e^{\pi\imath z(z-Q)}.
\end{equation}

Let us also introduce a closely related function $g_{b}(z)$, eq.(3.9) in \cite{BT} by the formula:

\begin{equation}\label{g_b in terms of G_b}
g_{b}(z) = \frac{\bar{\zeta}_{b}}{G_{b}(\frac{Q}{2}+\frac{1}{2\pi\imath b}\log z)}.
\end{equation}

Let $A$, $B$ be a pair of self-adjoint operators with the commutation relation $[A,B] = 2\pi\imath b^{2}$ and let $q = e^{\pi\imath b^{2}}$. Then one can define the following positive operators $u = e^{A}$, $v = e^{B}$ with the following commutation relation
\begin{equation}
    uv = q^{2}vu.
\end{equation}
For such operators the following identity ($q$-binomial theorem, see Appendix B in \cite{BT}) holds
\begin{equation}\label{q-binomial-1}
(u+v)^{\imath s} = \int\limits_{\mathcal{C}} d\tau \frac{G_{b}(-\imath b\tau)G_{b}(-\imath bs+\imath b\tau)}{G_{b}(-\imath bs)}u^{\imath s-\imath\tau}v^{\imath\tau},
\end{equation}
where the contour $\mathcal{C}$ goes along the real axis above the sequences of poles going down and below sequences of poles going up.

\section{Explicit realization of algebra of arbitrary complex devided powers of $U_{q}(\mathfrak{sl}(2))$}

Let $H$, $(K = q^{H})$, $E$, $F$ be generators of $U_{q}(\mathfrak{sl}(2))$ and let us introduce the following rescaled generators:
\begin{equation}\label{rescaled E}
    \mathcal{E} = -\imath (q-q^{-1})E,
\end{equation}
\begin{equation}\label{rescaled F}
    \mathcal{F} = -\imath (q-q^{-1})F,
\end{equation}
which satisfy the following relations
\begin{equation}
K\mathcal{E} = q^{2}\mathcal{E}K,
\end{equation}
\begin{equation}
K\mathcal{F} = q^{-2}\mathcal{F}K,
\end{equation}
\begin{equation}
\mathcal{E}\mathcal{F} - \mathcal{F}\mathcal{E} = -(q-q^{-1})(K-K^{-1}).
\end{equation}

\begin{prop}
The following operators define a representation of $U_{q}(\mathfrak{sl}(2))$
\begin{equation}
    H = 2\imath b^{-1}u,
\end{equation}
\begin{equation}\label{representation K}
    K = e^{-2\pi bu},
\end{equation}
\begin{equation}\label{representation E}
    \mathcal{E}  = e^{-\pi b(u-\alpha) + \imath b\partial_{u}} + e^{\pi b(u-\alpha) + \imath b\partial_{u}},
\end{equation}
\begin{equation}\label{representation F}
    \mathcal{F}  = e^{\pi b(u+\alpha) - \imath b\partial_{u}} + e^{-\pi b(u+\alpha) - \imath b\partial_{u}},
\end{equation}
\end{prop}
Note, that operators $\mathcal{E}$ and $\mathcal{F}$ have the form
\begin{equation}
\mathcal{E} =  U_{1} + V_{1},
\end{equation}
\begin{equation}
\mathcal{F} = U_{2} + V_{2},
\end{equation}
where $U_{1} = e^{-\pi b(u-\alpha) + \imath b\partial_{u}}$, $V_{1} = e^{\pi b(u-\alpha) + \imath b\partial_{u}}$,
$U_{2} = e^{\pi b(u+\alpha) - \imath b\partial_{u}}$, $V_{2} = e^{-\pi b(u+\alpha) - \imath b\partial_{u}}$, and
$U_{j}V_{j} = q^{2}V_{j}U_{j}$, $j = 1,2$. This allows one to define their arbitrary powers by (\ref{q-binomial-1}).
Let us introduce the notion of divided powers of $\mathcal{E}$ and $\mathcal{F}$ by the formulas
\begin{equation}\label{complex devided power}
    \mathcal{E}^{(\imath s)} = G_{b}(-\imath bs)\mathcal{E}^{\imath s},
\end{equation}
\begin{equation}\label{complex devided power F}
    \mathcal{F}^{(\imath t)} = G_{b}(-\imath bt)\mathcal{F}^{\imath t}.
\end{equation}
Then we have

\begin{prop}(Lemma 4.3, \cite{Ip}). Let the rescaled generators $K$, $\mathcal{E}$, $\mathcal{F}$ of $U_{q}(\mathfrak{sl}(2))$ be realized by the formulas (\ref{representation K})-(\ref{representation F}). Then for their powers we have the following realization
\begin{equation}\label{representation K^{ip}}
    K^{\imath p} = e^{-2\pi\imath bpu},
\end{equation}
\begin{equation}\label{representation E^{is}}
    \mathcal{E}^{(\imath s)} =
    e^{\frac{\pi\imath b^{2}s^{2}}{2} - \pi\imath bs(u-\alpha)}
    \frac{G_{b}(-\imath bs)G_{b}(\frac{Q}{2} + \imath\alpha +\imath bs -\imath u)}{G_{b}(\frac{Q}{2} + \imath\alpha -\imath u)} e^{-bs\partial_{u}},
\end{equation}
\begin{equation}\label{representation F^{it}}
    \mathcal{F}^{(\imath t)} = e^{\frac{\pi\imath b^{2}t^{2}}{2} + \pi\imath bt(u+\alpha)}
    \frac{G_{b}(-\imath bt)G_{b}(\frac{Q}{2} + \imath\alpha + \imath bt+\imath u)}{G_{b}(\frac{Q}{2}+\imath\alpha +\imath u)}e^{bt\partial_{u}}.
\end{equation}
\end{prop}
$\proof$
See the proof in \cite{Ip}, or for another proof one can use the $q$-binomial Theorem (\ref{q-binomial-1}) and after that apply $\tau$-binomial integral (\ref{tau-integral}).
$\Box$

We are going to consider the algebra ${A}(\mathfrak{sl}(2))$ spanned by the elements $\mathcal{E}^{(\imath s)}$, $\mathcal{F}^{(\imath t)}$, $K^{\imath p}$.

% \section{Hopf algebra of arbitrary complex devided powers of $U_{q}(\mathfrak{sl}(2))$ and its realization}

\begin{te}
For the operators $K^{\imath p}$, $\mathcal{E}^{(\imath s)}$, $\mathcal{F}^{(\imath t)}$ given by the formulas (\ref{representation K^{ip}})-(\ref{representation F^{it}}) the following commutation relations hold
\begin{equation}\label{KK relation}
    K^{\imath p_{1}}K^{\imath p_{2}} = K^{\imath p_{2}}K^{\imath p_{1}} = K^{\imath p_{1}+\imath p_{2}},
\end{equation}
\begin{equation}   \label{EE relation}
\mathcal{E}^{(\imath s_{1})}\mathcal{E}^{(\imath s_{2})} = \frac{G_{b}(-\imath bs_{1})G_{b}(-\imath bs_{2})}{G_{b}(-\imath bs_{1}-\imath bs_{2})}\mathcal{E}^{(\imath s_{1}+\imath s_{2})} = \mathcal{E}^{(\imath s_{2})}\mathcal{E}^{(\imath s_{1})},
\end{equation}
\begin{equation}\label{FF relation}
\mathcal{F}^{(\imath t_{1})}\mathcal{F}^{(\imath t_{2})} = \frac{G_{b}(-\imath bt_{1})G_{b}(-\imath bt_{2})}{G_{b}(-\imath bt_{1}-\imath bt_{2})}\mathcal{F}^{(\imath t_{1}+\imath t_{2})} = \mathcal{F}^{(\imath t_{2})}\mathcal{F}^{(\imath t_{1})},
\end{equation}
\begin{equation}    \label{K E relation}
K^{\imath p}\mathcal{E}^{(\imath s)} = e^{-2\pi\imath b^{2}ps}\mathcal{E}^{(\imath s)}K^{\imath p},
\end{equation}
\begin{equation}          \label{K F relation}
K^{\imath p}\mathcal{F}^{(\imath t)} = e^{2\pi\imath b^{2}pt}\mathcal{F}^{(\imath t)}K^{\imath p},
\end{equation}

\begin{equation}\label{generalized Kac's identity}
\mathcal{E}^{(\imath s)}\mathcal{F}^{(\imath t)} = \int\limits_{\mathcal{C}} d\tau e^{\pi bQ\tau}\mathcal{F}^{(\imath t+\imath\tau)}K^{-\imath \tau}
\frac{G_{b}(\imath b\tau)G_{b}(-bH + \imath b(s+t+\tau))}{G_{b}(-bH+\imath b(s+t+2\tau))}\mathcal{E}^{(\imath s + \imath \tau)},
\end{equation}
where the contour $\mathcal{C}$ goes along the real axis above the sequences of poles going down and below the sequences of poles going up.
\end{te}
The relations (\ref{KK relation}) are trivial.
Commutation relations (\ref{K E relation}), (\ref{K F relation}) trivially follow from the relation $e^{au}e^{b\partial_{u}} = e^{-ab}e^{b\partial_{u}}e^{au}$.
Relations (\ref{EE relation}), (\ref{FF relation}), (\ref{generalized Kac's identity}) are checked in the following two lemmas.

\begin{lem}
Let generators $\mathcal{E}^{(\imath s)}$, $\mathcal{F}^{(\imath t)}$ be realized by the operators from the previous proposition. Then these operators satisfy commutation relations (\ref{EE relation}) and (\ref{FF relation}) .
\end{lem}
$\proof$

Let us check the relation $\mathcal{E}^{(\imath s_{1})}\mathcal{E}^{(\imath s_{2})} = \frac{G_{b}(-\imath bs_{1})G_{b}(-\imath bs_{2})}{G_{b}(-\imath bs_{1}-\imath bs_{2})}\mathcal{E}^{(\imath s_{1}+\imath s_{2})}$. Substituting the formula (\ref{representation E^{is}}) into the left hand side we obtain
$$
\mathcal{E}^{(\imath s_{1})}\mathcal{E}^{(\imath s_{2})} =
e^{\frac{\pi\imath b^{2}s_{1}^{2}}{2} - \pi\imath bs_{1}(u-\alpha)}
\frac{G_{b}(-\imath bs_{1})G_{b}(\frac{Q}{2} + \imath\alpha +\imath bs_{1} -\imath u)}{G_{b}(\frac{Q}{2} + \imath\alpha -\imath u)} e^{-bs_{1}\partial_{u}}\times
$$
$$
e^{\frac{\pi\imath b^{2}s_{2}^{2}}{2} - \pi\imath bs_{2}(u-\alpha)}
\frac{G_{b}(-\imath bs_{2})G_{b}(\frac{Q}{2} + \imath\alpha +\imath bs_{2} -\imath u)}{G_{b}(\frac{Q}{2} + \imath\alpha -\imath u)} e^{-bs_{2}\partial_{u}} =
$$
$$
e^{\frac{\pi\imath b^{2}s_{1}^{2}}{2} - \pi\imath bs_{1}(u-\alpha)}
\frac{G_{b}(-\imath bs_{1})G_{b}(\frac{Q}{2} + \imath\alpha +\imath bs_{1} -\imath u)}{G_{b}(\frac{Q}{2} + \imath\alpha -\imath u)}\times
$$
$$
e^{\frac{\pi\imath b^{2}s_{2}^{2}}{2} - \pi\imath bs_{2}(u-bs_{1}-\alpha)}
\frac{G_{b}(-\imath bs_{2})G_{b}(\frac{Q}{2} + \imath\alpha +\imath bs_{1} +\imath bs_{2} -\imath u)}{G_{b}(\frac{Q}{2} + \imath\alpha +\imath bs_{1}-\imath u)} e^{-b(s_{1}+s_{2})\partial_{u}} =
$$
$$
e^{\frac{\pi\imath b^{2}(s_{1}+s_{2})^{2}}{2}-\pi\imath b(s_{1}+s_{2})(u-\alpha)}
\frac{G_{b}(-\imath bs_{1})G_{b}(-\imath bs_{2})G_{b}(\frac{Q}{2}+\imath\alpha+\imath bs_{1}+\imath bs_{2}-\imath u)}{G_{b}(\frac{Q}{2}+\imath\alpha-\imath u)}e^{-b(s_{1}+s_{2})\partial_{u}} =
$$
$$
\frac{G_{b}(-\imath bs_{1})G_{b}(-\imath bs_{2})}{G_{b}(-\imath bs_{1}-\imath bs_{2})}
e^{\frac{\pi\imath b^{2}(s_{1}+s_{2})^{2}}{2}-\pi\imath b(s_{1}+s_{2})(u-\alpha)}\times
$$
$$
\times\frac{G_{b}(-\imath bs_{1}-\imath bs_{2})G_{b}(\frac{Q}{2}+\imath\alpha+\imath bs_{1}+\imath bs_{2}-\imath u)}{G_{b}(\frac{Q}{2}+\imath\alpha-\imath u)}e^{-b(s_{1}+s_{2})\partial_{u}} =
$$
$$
\frac{G_{b}(-\imath bs_{1})G_{b}(-\imath bs_{2})}{G_{b}(-\imath bs_{1}-\imath bs_{2})}
\mathcal{E}^{(\imath s_{1}+\imath s_{2})}.
$$
The relation
$$
\mathcal{F}^{(\imath t_{1})}\mathcal{F}^{(\imath t_{2})} = \frac{G_{b}(-\imath bt_{1})G_{b}(-\imath bt_{2})}{G_{b}(-\imath bt_{1}-\imath bt_{2})}\mathcal{F}^{(\imath t_{1}+\imath t_{2})},
$$
is checked analogously.
$\Box$

\begin{lem}
Let generators $K^{\imath p}$, $\mathcal{E}^{(\imath s)}$, $\mathcal{F}^{(\imath t)}$ be realized by the operators (\ref{representation K^{ip}})-(\ref{representation F^{it}}). Then for such operators the generalized Kac's identity (\ref{generalized Kac's identity}) holds.
\end{lem}
$\proof$
Substituting into the left hand side of generalized Kac's identity
$$
\mathcal{E}^{(\imath s)}\mathcal{F}^{(\imath t)} = \int\limits_{\mathcal{C}} d\tau e^{\pi bQ\tau}\mathcal{F}^{(\imath t+\imath\tau)}K^{-\imath \tau}
\frac{G_{b}(\imath b\tau)G_{b}(-bH + \imath b(s+t+\tau))}{G_{b}(-bH+\imath b(s+t+2\tau))}\mathcal{E}^{(\imath s + \imath \tau)},
$$
explicit realization of generators (\ref{representation K^{ip}})-(\ref{representation F^{it}}), we have
$$
\mathcal{E}^{(\imath s)}\mathcal{F}^{(\imath t)} =
$$
$$
e^{\frac{\pi\imath b^{2}s^{2}}{2} - \pi\imath bs(u-\alpha)}
\frac{G_{b}(-\imath bs)G_{b}(\frac{Q}{2} + \imath\alpha +\imath bs -\imath u)}{G_{b}(\frac{Q}{2} + \imath\alpha -\imath u)} e^{-bs\partial_{u}}
e^{\frac{\pi\imath b^{2}t^{2}}{2} + \pi\imath bt(u+\alpha)}\times
$$
$$
\frac{G_{b}(-\imath bt)G_{b}(\frac{Q}{2} + \imath\alpha + \imath bt+\imath u)}{G_{b}(\frac{Q}{2}+\imath\alpha +\imath u)}e^{bt\partial_{u}} =
e^{\frac{\pi\imath b^{2}s^{2}}{2} + \frac{\pi\imath b^{2}t^{2}}{2} - \pi\imath b^{2}st + \pi\imath b(t-s)u + \pi\imath b(s+t)\alpha}
\times
$$
$$
\frac{G_{b}(-\imath bs)G_{b}(-\imath bt)G_{b}(\frac{Q}{2} + \imath\alpha + \imath bs-\imath u)G_{b}(\frac{Q}{2} + \imath\alpha + \imath bt-\imath bs+\imath u)}
{G_{b}(\frac{Q}{2}+\imath\alpha-\imath u)G_{b}(\frac{Q}{2} + \imath\alpha -\imath bs+\imath u)}
e^{b(t-s)\partial_{u}}.
$$
The right hand side is given by:
$$
\int d\tau e^{\pi bQ\tau}\mathcal{F}^{(\imath t+\imath\tau)}K^{-\imath \tau}
\frac{G_{b}(\imath b\tau)G_{b}(-bH + \imath b(s+t+\tau))}{G_{b}(-bH+\imath b(s+t+2\tau))}\mathcal{E}^{(\imath s + \imath \tau)} =
$$
$$
\int d\tau e^{\pi bQ\tau} e^{\frac{\pi\imath b^{2}}{2}(t^{2}+\tau^{2}+2t\tau)}e^{\pi\imath b(t+\tau)(u+\alpha)}
\frac{G_{b}(-\imath bt-\imath b\tau)G_{b}(\frac{Q}{2} + \imath\alpha + \imath bt+\imath u+\imath b\tau)}{G_{b}(\frac{Q}{2} +\imath\alpha +\imath u)}
e^{b(t+\tau)\partial_{u}}\times
$$
$$
e^{2\pi\imath b\tau u}
\frac{G_{b}(\imath b\tau)G_{b}(-2\imath u +\imath bs+\imath bt+\imath b\tau)}{G_{b}(-2\imath u+\imath bs+\imath bt+2\imath b\tau)}
e^{\frac{\pi\imath b^{2}}{2}(s^{2}+\tau^{2}+2s\tau)} e^{-\pi\imath b(s+\tau)(u-\alpha)}\times
$$
$$
\frac{G_{b}(-\imath bs-\imath b\tau)G_{b}(\frac{Q}{2}+\imath\alpha +\imath bs-\imath u+\imath b\tau)}{G_{b}(\frac{Q}{2}+\imath\alpha-\imath u)}
e^{-b(s+\tau)\partial_{u}} =
$$
$$
e^{\frac{\pi\imath b^{2}s^{2}}{2} + \frac{\pi\imath b^{2}t^{2}}{2} -\pi\imath b^{2}st +\pi\imath b(t-s)u + \pi\imath b(s+t)\alpha}
\frac{G_{b}(\frac{Q}{2} +\imath\alpha +\imath bs-\imath bt-\imath u)}{G_{b}(\frac{Q}{2} +\imath\alpha+\imath u)G_{b}(-2\imath u +\imath bs-\imath bt)} \times
$$
$$
\Big\{\int d\tau e^{2\pi\imath b^{2}\tau^{2} + 2\pi b(\frac{Q}{2} +\imath\alpha+\imath bt+\imath u)\tau}
G_{b}(-\imath bs-\imath b\tau)G_{b}(-\imath bt-\imath b\tau)G_{b}(-2\imath u+\imath bs-\imath bt-\imath b\tau)\times
$$
$$
\frac{G_{b}(\frac{Q}{2} +\imath\alpha+\imath bt+\imath u+\imath b\tau)G_{b}(\imath b\tau)}{G_{b}(\frac{Q}{2}+\imath\alpha -\imath bt-\imath u-\imath b\tau)}\Big\}
e^{b(t-s)\partial_{u}} = I
$$
This integral converges as $\tau\rightarrow \pm\infty$ for any parameters and the contour is deformed if necessary to separate sequences of poles going up and sequences of poles going down. In this case one can apply $6-9$ identity (\ref{6-9}) with
$$
A = -\imath bs,
$$
$$
B = -\imath bt,
$$
$$
C = -2\imath u +\imath bs-\imath bt,
$$
$$
D = \frac{Q}{2} + \imath\alpha +\imath bt+\imath u,
$$
$$
A+B+C+D = \frac{Q}{2} +\imath\alpha -\imath bt -\imath u,
$$
to obtain
$$
I =
e^{\frac{\pi\imath b^{2}s^{2}}{2} + \frac{\pi\imath b^{2}t^{2}}{2} -\pi\imath b^{2}st +\pi\imath b(t-s)u + \pi\imath b(s+t)\alpha}
\frac{G_{b}(\frac{Q}{2} +\imath\alpha +\imath bs-\imath bt-\imath u)}{G_{b}(\frac{Q}{2} +\imath\alpha+\imath u)G_{b}(-2\imath u +\imath bs-\imath bt)}\times
$$
$$
G_{b}(-\imath bs)G_{b}(-\imath bt)G_{b}(-2\imath u+\imath bs-\imath bt)
\frac{G_{b}(\frac{Q}{2} +\imath\alpha+\imath bt-\imath bs+\imath u)G_{b}(\frac{Q}{2}+\imath\alpha+\imath u)G_{b}(\frac{Q}{2} +\imath\alpha+\imath bs-\imath u)}
{G_{b}(\frac{Q}{2}+\imath\alpha-\imath bs+\imath u)G_{b}(\frac{Q}{2}+\imath\alpha-\imath u)G_{b}(\frac{Q}{2}+\imath\alpha+\imath bs-\imath bt-\imath u)}
$$
$$
\times e^{b(t-s)\partial_{u}} =
e^{\frac{\pi\imath b^{2}s^{2}}{2} + \frac{\pi\imath b^{2}t^{2}}{2} - \pi\imath b^{2}st + \pi\imath b(t-s)u + \pi\imath b(s+t)\alpha}
\times
$$
$$
\frac{G_{b}(-\imath bs)G_{b}(-\imath bt)G_{b}(\frac{Q}{2} + \imath\alpha + \imath bs-\imath u)G_{b}(\frac{Q}{2} + \imath\alpha + \imath bt-\imath bs+\imath u)}
{G_{b}(\frac{Q}{2}+\imath\alpha-\imath u)G_{b}(\frac{Q}{2} + \imath\alpha -\imath bs+\imath u)}
e^{b(t-s)\partial_{u}} =
\mathcal{E}^{(\imath s)}\mathcal{F}^{(\imath t)}.
$$

$\Box$

We have explicitly checked all the commutation relations in the Hopf algebra, generated by complex powers of generators of $U_{q}(\mathfrak{sl}(2))$ in the positive principal series representations.

\section{Appendix}

\subsection{Quantum dilogarithm and its properties}
The basic properties of non-compact quantum dilogarithm/double sine listed below are extracted mainly from \cite{KLSTS}, \cite{BT}, \cite{V}.
Introduce the following notation $q = e^{\pi\imath b^{2}}$, $\tilde{q} = e^{\pi\imath b^{-2}}$, $Q = b+b^{-1}$, $\zeta_{b} = e^{\frac{\pi\imath}{4} + \frac{\pi\imath(b^{2}+b^{-2})}{12}}$.\\*
\\*
\textbf{The integral representation of $G_{b}(z)$:}
\begin{equation}
\log G_{b}(z) = \log\bar{\zeta}_{b} - \int\limits_{\mathbb{R}+\imath 0} \frac{dt}{t}\frac{e^{zt}}{(1-e^{bt})(1-e^{b^{-1}t})}.
\end{equation}
\textbf{Noncompact analog of $q$-exponential $g_{b}(z)$:}
\begin{equation}
g_{b}(z) = \frac{\bar{\zeta}_{b}}{G_{b}(\frac{Q}{2}+\frac{1}{2\pi\imath b}\log z)}.
\end{equation}
\textbf{Product representation:}
\begin{equation}
G_{b}(x) = \bar{\zeta}_{b}\frac{\prod\limits_{n=1}^{\infty}(1-e^{2\pi\imath b^{-1}(x-nb^{-1})})}{\prod\limits_{n=0}^{\infty}(1-e^{2\pi\imath b(x+nb)})},
\end{equation}
\begin{equation}
g_{b}(x) = \frac{\prod\limits_{n=0}^{\infty}(1+xq^{2n+1})}{\prod\limits_{n=0}^{\infty}(1+x^{b^{-2}}\tilde{q}^{-2n-1})}.
\end{equation}
\textbf{Functional equations:}
\begin{equation}
G_{b}(x +b^{\pm 1}) = (1-e^{2\pi\imath b^{\pm 1}x})G_{b}(x),
\end{equation}
or more generally
\begin{equation} \label{func eq}
\frac{G_{b}(x+n_{1}b+n_{2}b^{-1})}{G_{b}(x)} =
\prod\limits_{k_{1} =0}^{n_{1}-1}(1-q^{2k_{1}}e^{2\pi\imath bx})\prod\limits_{k_{2} =0}^{n_{2}-1}(1-\tilde{q}^{2k_{2}}e^{2\pi\imath b^{-1}x}),
\end{equation}
\begin{equation}
g_{b}(q^{-1}x) = (1+x)g_{b}(qx).
\end{equation}
\textbf{Reflection formula:}
\begin{equation}\label{reflection}
G_{b}(x)G_{b}(Q-x) = e^{\pi\imath x(x-Q)}.
\end{equation}
\textbf{Poles and zeros:}
\begin{equation}\begin{split}\label{Poles}
\lim_{x\rightarrow 0}xG_{b}(x-n_{1}b-n_{2}b^{-1}) = \frac{1}{2\pi}\prod\limits_{k_{1}=1}^{n_{1}}(1-q^{-2k_{1}})^{-1}\prod\limits_{k_{2}=1}^{n_{2}}(1-\tilde{q}^{-2k_{2}})^{-1},    \\
\lim_{x\rightarrow 0}xG_{b}^{-1}(x+Q+n_{1}b+n_{2}b^{-1}) =      \\
\frac{1}{2\pi}(-1)^{n_{1}+n_{2}+1}q^{-n_{1}(n_{1}+1)}\tilde{q}^{-n_{2}(n_{2}+1)}\prod\limits_{k_{1}=1}^{n_{1}}(1-q^{-2k_{1}})^{-1}\prod\limits_{k_{2}=1}^{n_{2}}(1-\tilde{q}^{-2k_{2}})^{-1}.
\end{split}\end{equation}
\textbf{Tau-binomial integral} \cite{FKV},\cite{Ka},\cite{PT}:
\begin{equation}\label{tau-integral}
\int\limits_{\mathcal{C}} d\tau e^{-2\pi b\beta\tau}\frac{G_{b}(\alpha+\imath b\tau)}{G_{b}(Q+\imath b\tau)} =
\frac{G_{b}(\alpha)G_{b}(\beta)}{G_{b}(\alpha+\beta)},
\end{equation}
where the contour $\mathcal{C}$ goes along the real axis above the sequences of poles going down and below sequences of poles going up.\\*
\\*
\textbf{6-9 identity \cite{V}:}
\begin{equation}\begin{split}        \label{6-9}
\frac{G_{b}(A)G_{b}(B)G_{b}(C)G_{b}(A+D)G_{b}(B+D)G_{b}(C+D)}
{G_{b}(A+B+D)G_{b}(A+C+D)G_{b}(B+C+D)} = \\
\int\limits_{\mathcal{C}} d\tau e^{2\pi\imath\tau^{2}-2\pi D\tau}
\frac{G_{b}(A+\imath\tau)G_{b}(B+\imath\tau)G_{b}(C+\imath\tau)G_{b}(D-\imath\tau)G_{b}(-\imath\tau)}{G_{b}(A+B+C+D+\imath\tau)},
\end{split}\end{equation}
where the contour $\mathcal{C}$ goes along the real axis above the sequences of poles going down and below sequences of poles going up.\\*
\\*
\textbf{$q$-binomial theorem \cite{BT}:}\\*
Let $u$, $v$ be positive self-adjoint operators subject to the relations $uv = q^{2}vu$.
Then:
\begin{equation}\label{q-binomial}
(u+v)^{\imath s} = \int\limits_{\mathcal{C}} d\tau \frac{G_{b}(-\imath b\tau)G_{b}(-\imath bs+\imath b\tau)}{G_{b}(-\imath bs)}u^{\imath s-\imath\tau}v^{\imath\tau},
\end{equation}
where the contour $\mathcal{C}$ goes along the real axis above the sequences of poles going down and below sequences of poles going up.

\end{document}